\documentclass[preprint,11pt]{elsarticle}
\usepackage{graphicx}
\usepackage{amsfonts}
\usepackage{amsmath}

\usepackage{amscd,amsthm}

\usepackage{alltt}
\usepackage{array}
\usepackage{multirow}
\usepackage{lscape}
\newcommand{\dexp}[1]{\times 10^{#1}}

\def\H{\mathcal{H}}

\newtheorem{thm}{Theorem}

\newtheorem{rmk}[thm]{Remark}
\journal{Physica D}

\begin{document}

\begin{frontmatter}

\title{Equilibrium search algorithm of a perturbed quasi-integrable system}
%\subtitle{\texttt{N}amur \texttt{A}lgorithm \texttt{F}or \texttt{F}orced \texttt{O}scillations}


%\author{Beno\^it Noyelles         \and
%        Nicolas Delsate \and
%	Timoteo Carletti
%}


%\institute{Namur Center for Complex Systems, naXys, University of Namur\\
%              Rempart de la Vierge 8 -- B-5000 Namur -- BELGIUM\\
%              Tel.: +3281724940\\
%              \email{benoit.noyelles@fundp.ac.be}             \\
%   B. Noyelles is F.R.S.-FNRS Post-doctoral research fellow
%}


%\maketitle

\author{Beno\^it Noyelles\fnref{fn1}}
\ead{benoit.noyelles@fundp.ac.be}
\author{Nicolas Delsate}
\ead{nicolas.delsate@math.fundp.ac.be}
\author{Timoteo Carletti}
\ead{timoteo.carletti@fundp.ac.be}
\address{University of Namur -- Department of Mathematics \& NAXYS \\
    Rempart de la Vierge 8 -- B-5000 Namur -- BELGIUM}
\fntext[fn1]{F.R.S.-FNRS post-doctoral research fellow}

\begin{abstract}
We hereby introduce and study an algorithm able to search for initial conditions corresponding to orbits presenting forced oscillations terms only, namely to completely remove the free or proper oscillating part.

This algorithm is based on the Numerical Analysis of the Fundamental Frequencies algorithm by J. Laskar, for the identification of the free and forced oscillations, the former being iteratively removed from the solution by carefully choosing the initial conditions.

We proved the convergence of the algorithm under suitable assumptions, satisfied in the Hamiltonian framework whenever the d'Alembert characteristic holds true. In this case, 
with polar canonical variables, we also proved that this algorithm converges quadratically. We provided a relevant application: the forced prey-predator problem.

\end{abstract}


\end{frontmatter}

\section{Introduction}

\par A large number of low-dissipative problems can be conveniently modelled over a quite short time interval as if they were conservative systems, 
i.e. in neglecting the short term influence of the dissipative processes, and moreover to evolve close to a stable equilibrium because the oscillations 
around it can be assumed to have been damped in the past. This is often the case in celestial mechanics, where for instance the resonant rotation 
of planetary satellites is assumed to be at an equilibrium, and thus only depending on an external torque, without any free oscillation. However, 
the problem is more generally applicable to a large class of quasiperiodically forced systems, for instance we will provide an application to a 
well-known prey-predator system subjected to an external periodic forcing, accounting for the seasonal changes. Let us observe that the method is 
far more general and the applicability domains goes well beyond the two examples hereby presented as it can be seen.

More precisely we consider systems having stable equilibria and quasi-periodic orbits in their neighborhoods. We start from a system defined by an ordinary differential equation alike

\begin{equation}
  \dot{\vec{X}}(t)=f(\vec{X})+g(\vec{X},t),
\end{equation}
where $f:\mathbb{B}\subset \mathbb{C}^n \rightarrow \mathbb{C}^n$, $g:\mathbb{B}\subset \mathbb{C}^n \times \mathbb{R} \rightarrow \mathbb{C}^n$ being an external perturbation, and $\vec{X}\in \mathbb{B}\subset \mathbb{C}^n$, $\mathbb{B}$ being an open subset of $\mathbb{C}^n$ containing the equilibrium we are interested in, and we assume that every component $x_i(t)$ of the solution, $\vec{X}(t)$, in a neighborhood of the equilibrium can be expressed as a convergent infinite sum of periodic trigonometric monomials, such as: 
\begin{equation}
  x_i(t)=\sum_l a_l\exp\left(\imath f_lt\right)\, ,
  \label{eq:x}
\end{equation}
where $a_l$ are complex amplitudes and $f_l$ are real frequencies. We can split these frequencies into 2 groups: free, or proper frequencies of the \lq\lq unperturbed\rq\rq system, i.e. when $g\equiv 0$, and forced frequencies due to the external perturbation $g$.

This decomposition holds for quasi-integrable systems in the Hamiltonian framework thanks to the Kolmogorov-Arnold-Moser (KAM) theorem~\cite{a63b,k54,m62}, that ensures the existence of quasi-periodic invariant tori under suitable hypotheses: the non--degeneracy of the integrable part of the Hamiltonian function, the Diophantine condition on the frequencies and the smallness of the perturbation; let us observe that our setting is a simplified version of the P\"oschel result~\cite{poschel1989} where the normal frequencies, hereby called proper frequencies, do not depend on the torus frequencies, hereby named forced ones. On the other hand, the theorem by Nekhoroshev~\cite{n77,n79} proves that the orbits of a slightly perturbed system stay close to the orbits of the unperturbed system over an exponentially long time with respect to the inverse of the perturbation, and thus~\eqref{eq:x} is a good approximation for any realistic physical times.

For example, in the context of the resonant rotation of the Moon around the Earth, we can split the frequencies $f_l$ into two groups: the frequencies $\nu_j$ that are forced perturbations displacing the equilibrium, and the frequencies $\omega_i$ that are free librations of the system around the equilibrium. The frequencies $\nu_j$ are due to an external forcing (in our example, the gravitational perturbations of the Sun and the planets on the rotation of the Moon), while $\omega_i$ depend on the intrinsic properties of the system (in our example, on the gravity field of the Moon). Slow dissipative processes acting since the origin of the system (4.5 Gyr in the present case) are expected to have damped the free librations (see e.g. \cite{p05} for the rotation of Mercury, and \cite{gdnr96} for the planetary satellites). 

The numerical computation of the wanted orbits, namely whose free librations have negligible amplitudes, i.e. the most realistic ones, requires the accurate knowledge of 
the initial conditions corresponding to the equilibrium. This task can be tricky whenever the system is subject to several different sinusoidal perturbations. In such a case, an analytical determination of the equilibrium would either give a too approximate solution because related to a too simple system, or would be too difficult to use for the complete system to obtain the required accuracy.

A straightforward way to overcome this difficulty is to use a numerical algorithm that iteratively corrects the initial conditions by identifying and then
removing the free oscillations around the equilibrium. The idea of using the frequency analysis to refine the initial conditions is not completely new in the literature, in fact it has already been used by Locatelli \& Giorgili \cite{lg00} for a computer-assisted proof of the KAM theorem. Here, in removing the free oscillations, we aim at reducing the numbers of dimensions of the solutions, i.e. we are looking for solutions lying on lower dimension tori with respect to the full dimension $n+p$. This has already been done by Noyelles et al. (see e.g. \cite{dnrl09,n09,ndl10}) and Robutel et al. \cite{rrc11} for representing the resonant spin-orbit rotation 3:2 or 1:1 in the Solar System, by Couetdic et al.~\cite{c09} in the framework of (exo)planetary systems, and by Delsate \cite{d11} for the dynamics of a spacecraft around the asteroid Vesta. The aim of our paper is to clearly state the algorithm, discuss its applicability domains and moreover prove under suitable hypotheses its convergence.

The paper is organized as follows. Section~\ref{sec:naffo} is devoted to the presentation of the algorithm and the numerical tools needed. Then 
in Section~\ref{sec:proof} we will present a convergence proof of the method in the Hamiltonian framework, before addressing its extension
to a more general framework in Section~\ref{sec:generalcase}. In Section~\ref{sec:applbio} an application of the algorithm will be presented, to 
a problem from mathematical biology. We finally review briefly other algorithms in Section~\ref{sec:alternate}, before summing up and drawing 
our conclusions in Section~\ref{sec:concl}.

\section{The algorithm}
\label{sec:naffo}

\par The goal of this section is to present an algorithm to search for the equilibrium of a periodically forced system. Let observe that its aim is not to prove the existence of the wanted orbit that should be assumed to exist a priori, but to determine the initial condition providing the wanted orbit. For the conditions of existence of such orbits, we refer the interested reader to e.g. \cite{js96}.

The algorithm requires hence the existence of quasi-periodic orbits around the equilibrium and thus it relies on a method able to reconstruct these quasi-periodic orbits by identifying the forced and the free terms in the frequency space. To accomplish this task, in the following we decided to use the NAFF algorithm, Numerical Analysis of the Fundamental Frequencies, by Laskar~\cite{l93,l05}, shortly described in Appendix A. We used this algorithm because it has proved its efficiency, but any other algorithm allowing to detect accurately the frequencies could be used as well.

Let us consider a system described by the differential equation
\begin{equation}
  \label{eq:edo}
  \dot{\vec{X}}=f(\vec{X})\, ,
\end{equation}
in a neighborhood of an either (quasi)periodic or static solution with frequencies $(\omega_i)_{1\leq i\leq m}$ with $0\leq m\leq n$ written as $\vec{\omega}\in\mathbb{R}^n$, hereby named the equilibrium of the free motion, $\vec{X}$ being a n-dimensional vector. Let us then modify the system by adding a quasiperiodic forced term:
\begin{equation}
  \label{eq:edo2}
  \dot{\vec{X}}=f(\vec{X})+g(\vec{X},t)\, ,
\end{equation}
where $g(\vec{X},t)$ is quasiperiodic, its frequencies being $(\nu_j)_{1\leq j\leq p}$, also written as $\vec{\nu}\in\mathbb{R}^p$. So, the functions $f$ and $g$ are defined from an open set of $\mathbb{C}^n$ and $\mathbb{C}^n\times\mathbb{R}$ respectively, to $\mathbb{C}^n$, and $n$ and $p$ are strictly positive integers.

 Let us assume that $\omega_i$ and $\nu_j$ are rationally independent and thus secular terms are not allowed in the solution. Let us also denote the solution with a generic initial condition $\vec{X}$ by $\vec{\phi}(t;\vec{X}):\mathbb{R}\times\mathbb{C}^n\rightarrow\mathbb{C}^n$. Using NAFF we can identify the contribution of each frequency and thus obtain the decomposition
\begin{equation}
  \label{eq:sol}
  \vec{\phi}\left(t;\vec{X}\right)=\sum_{\vec{l}\in\mathbb{Z}^n,\vec{m}\in\mathbb{Z}^p}\vec{\phi}_{\vec{l}\vec{m}}\left(\vec{X}\right)e^{\imath (\vec{\omega}\cdot\vec{l} +\vec{\nu}\cdot\vec{m})t}\, ,
\end{equation}
that can be rewritten by separating free from forced oscillating terms, as follows

\begin{eqnarray}
  \vec{\phi}\left(t;\vec{X}\right) & = & \sum_{\vec{m}\in\mathbb{Z}^p}\vec{\phi}_{\vec{0},\vec{m}}\left(\vec{X}\right)e^{\imath\vec{\nu}\cdot\vec{m}t}+
\sum_{\vec{l}\ne \vec{0}, \vec{m}\in\mathbb{Z}^p}\vec{\phi}_{\vec{l},\vec{m}}\left(\vec{X}\right)e^{\imath(\vec{\omega}\cdot\vec{l}+\vec{\nu}\cdot\vec{m}) t} \nonumber \\
 & =: & \vec{S}(t;\vec{X})+\vec{L}(t;\vec{X})\, \label{eq:sol2} ,
\end{eqnarray}
where the quantities $S$ and $L$ have been here defined by identifying left and right hand sides, they are functions defined from an open set of 
$\mathbb{C}^n\times\mathbb{R}$ to $\mathbb{C}^n$. Their uniqueness comes from the assumption that the determination of the quasi-periodic series has a perfect accuracy. 

The assumption of the existence of a $\left(2\pi/\nu_j\right)_{1\leq j \leq p}$ quasiperiodic orbit, i.e. lying on a p-torus, namely composed only by 
forced oscillating terms, translates into the existence of an isolated initial condition $\vec{X}_{\infty}$ such that the solution 
with initial condition $\vec{X}_{\infty}$ is $\left(2\pi/\nu_j\right)_{1\leq j \leq p}$ quasiperiodic, namely:

\begin{equation*}
  \vec{\phi}\left(t;\vec{X}_{\infty}\right)=\vec{S}(t;\vec{X}_{\infty})\, ,
\end{equation*}
or equivalently $\vec{L}(t;\vec{X}_{\infty})\equiv \vec{0}$. $\vec{X}_{\infty}$ is unknown and our algorithm is a way to determine it. \\

This algorithm can be formulated as follows:

\begin{enumerate}

\item Take $\vec{X}_0$ sufficiently close to $\vec{X}_{\infty}$;

\item Integrate the ODE~\eqref{eq:edo2}, thus determine the solution
  $\vec{\phi}\left(t;\vec{X}_0\right)$ 
  and then, using NAFF, get the decomposition
  $\vec{\phi}\left(t;\vec{X}_0\right)=\vec{S}(t;\vec{X}_0)+\vec{L}(t;\vec{X}_0)$;  

\item Define $\vec{X}_1=\vec{S}(0;\vec{X}_0)$ and iterate point 2 using $\vec{X}_1$ instead of $\vec{X}_0$;

\item In this way the algorithm will produce a sequence $\vec{X}_n$, iteratively defined by  

\begin{equation}
  \label{eq:iteration}
  \vec{X}_{n+1}=\vec{S}(0;\vec{X}_n)\, ,   
\end{equation}
such that, when $n \rightarrow \infty$, $\vec{X}_n\rightarrow \vec{X}_{\infty}$, or equivalently $\vec{L}(t;\vec{X}_n)\rightarrow \vec{0}$. 
\end{enumerate}

\begin{rmk}
 Once numerically implemented, we can define at least two ways to define a stopping criterion for this algorithm: we can consider the convergence as reached when the largest amplitude associated with the free oscillations is too small to be detected, or when it is small enough to not significantly alter the determination of the forced oscillations. 
\end{rmk}

We here make no hypothesis on the dimensions $n$ and $p$. This algorithm has already been used in the following contexts:

\begin{itemize}
 
\item Noyelles (2009) \cite{n09}: $n=3$, $p=13$
\item Dufey et al. (2009) \cite{dnrl09}: $n=2$, $p=5$
\item Couetdic et al. (2010) \cite{c09}: $n=1$, $p=3$
\item Noyelles et al. (2010) \cite{ndl10}: $n=4$, $p=5$
\item Robutel et al. (2011) \cite{rrc11}: $n=1$, $p=6$
\item Delsate (2011) \cite{d11}: $n=1$, $p=3$
\end{itemize}

Its convergence under some assumptions will be proved in the next two sections.

\section{Convergence of the algorithm in the Hamiltonian framework}
\label{sec:proof}

%\par In this section, we tackle the problem of the convergence of the algorithm. We first prove it in the Hamiltonian framework, where the 
%d'Alembert characteristic \cite{h74} gives us also the convergence rate. Then we address the possible extensions to a wider framework

%\subsection{In the Hamiltonian framework}

Let us consider a $(n+p)$-DOF Hamiltonian ${H}(y_i,x_i,\Lambda_j,\lambda_j)$, $y_i$ ($1\leq i\leq n$) and $\Lambda_j$ 
($1\leq j\leq p$) being the momenta and $x_i$ and $\lambda_j$ the variables. We assume that the forcing is due to 
$\Lambda_j,\lambda_j$ and that they are respectively actions and angles variables. Moreover, we hypothesize that the Hamiltonian 
system can be locally described by a forced perturbed harmonic oscillator:

\begin{equation}
  \label{eq:birkhoff}
  \H(U_i,u_i,\Lambda_j,\lambda_j)=\sum_i\omega_i U_i + \epsilon \H_1(U_i,u_i,\Lambda_j,\lambda_j)\, ,
\end{equation}
where $\omega_i$ is the frequency of the free oscillations and $\epsilon \H_1$ is a small perturbation. This is a degenerate 
setting where existence of invariant tori can be studied thanks to Birkhoff theory (see e.g. \cite{m68}). The Hamiltonian $\H$ has 
been obtained after canonical transformations among which an untangling one that removes the second-order cross terms alike 
$y_jy_k$ with $j\ne k$ \cite{hl05}, and the classical canonical polar transformation: 

\begin{eqnarray}
%\begin{cases}
  x_i  &=& \sqrt{2U_iZ_i}\sin u_i, \nonumber \\
  y_i  &=&  \sqrt{2U_i/Z_i}\cos u_i, \label{eq:polar}
%\end{cases}
\end{eqnarray}
where $Z_i$ is a constant term, suitably chosen so that the first order terms in $\sqrt{U_i}$ disappear. The transformation 
(\ref{eq:polar}) can be seen as the expression of the variables $x_i$ and $y_i$ after averaging with respect to the forced 
perturbation, i.e.  

\begin{eqnarray}
\begin{cases}
  x_i(U_i,u_i,-,-) & =  \sqrt{2U_iZ_i}\sin u_i, \nonumber \\
  y_i(U_i,u_i,-,-) & =  \sqrt{2U_i/Z_i}\cos u_i\, . \label{eq:polar2}
\end{cases}
\end{eqnarray}

There are several perturbative methods, based on the hypotheses that the perturbation is small and far enough from resonances with 
the proper frequencies $\omega_i$, that allow to derive the complete expression of $x_i$ and $y_i$, i.e. including their 
dependencies on $\Lambda_j$ and $\lambda_j$. We hereby propose to use the method of the Lie transforms (see \cite{d69} or 
\cite{d10} for detailed explanations). In general, the formal series do not converge, nevertheless, such perturbation
techniques can be justified through Poincar\'e's theory of asymptotic approximations \cite{p99}. It is the difference
between the convergence seen by geometers, on a finite-size interval, and the one seen by astronomers, who see the series
as asymptotic expansions around $0$.

Under the assumption of an analytic Hamiltonian $\H$, Henrard \cite{h74} has proved that at the $N$--th step of the Lie transforms algorithm, the 
functions $x_{q,N}$ and $y_{q,N}$ with $1\leq q \leq n$, follow the d'Alembert characteristic for $(U_{i,N},u_{i,N})$, i.e. they can be written as 

\begin{equation}
  \label{eq:xfour}
  \begin{split}
  x_{q,N}(U_{i,N},u_{i,N},\Lambda_j,\lambda_j)=\sum_{k_i,h_j\in\mathbb{Z}}\alpha_{k_i,h_j}\prod_{k_i} U_{i,N}^{\frac{|k_i|}{2}}
  \Big(1+\sum_{m\geq1}\gamma_{m_i,k_i,h_j} \prod U_{i,N}^{m_i}\Big) \\
  \times \exp\Big(\imath\big(\sum_{k_i}k_iu_{i,N}+\sum_{h_j}\lambda_j\big)\Big)
  \end{split}
\end{equation}
and
\begin{equation}
  \label{eq:yfour}
  \begin{split}
  y_{q,N}(U_{i,N},u_{i,N},\Lambda_j,\lambda_j)=\sum_{k_i,h_j\in\mathbb{Z}}\beta_{k_i,h_j}\prod_{k_i} U_{i,N}^{\frac{|k_i|}{2}}
  \Big(1+\sum_{m\geq1}\delta_{m_i,k_i,h_j} \prod U_{i,N}^{m_i}\Big) \\
  \times \exp\Big(\imath\big(\sum_{k_i}k_iu_{i,N}+\sum_{h_j}\lambda_j\big)\Big),
  \end{split}
\end{equation}
where $\alpha_{k_i,h_j}$, $\beta_{k_i,h_j}$, $\gamma_{m_i,k_i,h_j}$ and $\delta_{m_i,k_i,h_j}$ are complex constants, 
and $m_i$, $k_i$ and $h_j$ are integers. The right-hand members are assumed to converge absolutely for $U_{i,N}\leq U_i^{\bullet}$ for some 
positive $U_i^{\bullet}$.

So we get the time evolution of the above functions given by:

\begin{equation}
  \begin{split}
  x_{q,N}(t)=\sum_{k_i,h_j\in\mathbb{Z}}\alpha_{k_i,h_j}\prod_{k_i} U_{i,N}(t)^{\frac{|k_i|}{2}}
  \Big(1+\sum_{m\geq1}\gamma_{m_i,k_i,h_j} \prod U_{i,N}(t)^{m_i}\Big) \\
  \times \exp\Big(\imath\big(\sum_{k_i}k_iu_{i,N}(t)+\sum_{h_j}h_j\lambda_j(t)\big)\Big)
  \end{split}
\end{equation}

\begin{equation}
  \begin{split}
  y_{q,N}(t)=\sum_{k_i,h_j\in\mathbb{Z}}\beta_{k_i,h_j}\prod_{k_i} U_{i,N}(t)^{\frac{|k_i|}{2}}
  \Big(1+\sum_{m\geq1}\delta_{m_i,k_i,h_j} \prod U_{i,N}(t)^{m_i}\Big) \\
  \times \exp\Big(\imath\big(\sum_{k_i}k_iu_{i,N}(t)+\sum_{h_j}h_j\lambda_j(t)\big)\Big).
  \end{split}
\end{equation}

Assuming that the truncation order $M$ is large enough such that the functions $U_{i,N}(t)\sim U_{i,N}^{\prime}$ is almost constant, as they 
would be in the limit $M\rightarrow \infty$, we can obtain the decomposition of $x_q^N(t)$, and similarly for $y_q^N(t)$, into forced
and free oscillating terms as follows

\begin{equation}
  \begin{split}
  x_q^N(t)=\underbrace{\sum_{k_i,h_j\in\mathbb{Z}}\alpha_{0,h_j}\Big(1+\sum_{m\geq1}\gamma_{m_i,0,h_j} \prod U_{i,N}'^{m_i}\Big)
  \exp\Big(\imath\big(\sum_{h_j}\lambda_j(t)\big)\Big)}_{S_{q,N}}+\\
\underbrace{\sum_{k_i\ne 0,h_j\in\mathbb{Z}}\alpha_{k_i,h_j}\prod_{k_i} U_{i,N}'^{\frac{|k_i|}{2}}
\Big(1+\sum_{m\geq1}\gamma_{m_i,k_i,h_j} \prod U_{i,N}'^{m_i}\Big)\exp\Big(\imath\big(\sum_{k_i}k_iu_{i,N}'+\sum_{h_j}\lambda_j(t)\big)\Big)}_{L_{q,N}}\, .
   \end{split}
   \label{eq:xn}
\end{equation}

For any $N$  the $S_{q,N}$ terms, components of $\vec{S}_{N}$ behave as $$S_{q,N}\sim A+\sum_{i,j}B_{i,j} \sqrt{U_{i,N}'}\sqrt{U_{j,N}'}+ \dots$$ 
while the $L_{q,N}$ terms, components of $\vec{L}_{q,N}$, behave as $$L_{q,N}\sim \sum_{i}C_i\sqrt{U_{i,N}'}+ \dots$$, where $"\ldots"$ means 
"higher order terms".

Let us observe that setting $t=0$ in the Eq.~\ref{eq:xn} we get from the very definition of our algorithm:
\begin{equation}
	x_{q,{N+1}} =  S_{q,N}(0)\, ,
	\label{eq:xq1}
\end{equation}
that should equals the left hand side of~\ref{eq:xn} written with $N+1$, that is:
\begin{equation}
	x_{q,{N+1}}  \sim  A+\sum_{i,j}B_{i,j} \sqrt{U_{i,N}'}\sqrt{U_{j,N}'}+ \dots \label{eq:xq2}
\end{equation}

\par From the Eq.\ref{eq:xq1} and \ref{eq:xq2} we get

\begin{equation}
	\label{eq:cvrate}
	\sqrt{U_{i,N+1}'}=\mathcal{O}(U_{i,N}).
\end{equation}
If we assume the different components $U_{i,N}$ of the vector $\vec{U_N}$ at the iteration $N$ to be of the same order of magnitude, then we have a quadratic convergence of the algorithm.

\section{The convergence in a more general framework}
\label{sec:generalcase}

\par In a general framework, the convergence cannot be checked a priori. We here investigate a condition leading to this convergence, before discussing its relevance.

\subsection{The hypothesis}

\par The convergence of the algorithm can be proved by showing the convergence
of the sequence $(\vec{X}_k)_{k\in\mathbb{N}}$ defined by~\eqref{eq:iteration}.

\par If we have $\vec{X}_k\rightarrow \vec{P}$ for some $\vec{P}$, $\vec{X}_k$ and $\vec{P}$ being of dimension $n$, then by 
rewriting~\eqref{eq:iteration} as follows
\begin{equation}
  \label{eq:iteration2}
  \vec{X}_{k+1}=S(0;\vec{X}_k)=\phi(0;\vec{X}_k)-L(0;\vec{X}_k)=\vec{X}_k-L(0;\vec{X}_k)\, ,   
\end{equation}
we straightforwardly get by continuity of $L$ and the convergence hypothesis of $\vec{X}_k$
\begin{equation*}
  \vec{0}=  \lim_{n\rightarrow \infty}\left(\vec{X}_{n+1}-\vec{X}_n\right)= -\lim_{n\rightarrow\infty}L(0;\vec{X}_n)=-L(0;\vec{P})\, ,
\end{equation*}
then by uniqueness of the $2\pi/\nu_j$-quasiperiodic orbit we can conclude that $\vec{P}=\vec{X}_{\infty}$.

We provide a proof of the convergence of the algorithm assuming the following hypothesis: the Jacobian matrices $\Sigma$ and $\Lambda$ 
of the functions $\vec{S}$ and $\vec{L}$ do satisfy 
\begin{equation}
  \label{eq:hypo}
   \lim_{\vec{X}\rightarrow \vec{X}_{\infty}}\left[\Lambda(\vec{X})\right]^{-1}\Sigma(\vec{X})=0\, .
\end{equation}

\subsection{Proof}

From the definition~\eqref{eq:iteration} of $\vec{X}_{k+1}$ and the very definition of $\vec{X}_{k+1}$ as initial condition of the 
orbit we get:

\begin{equation}
  \label{eq:sls}
  S(0;\vec{X}_k)= \vec{X}_{k+1}=\phi(0;\vec{X}_{k+1})=S(0;\vec{X}_{k+1})+L(0;\vec{X}_{k+1})\, .
\end{equation}
This relation defines implicitly the map that associates to $\vec{X}_k$ the next iteration $\vec{X}_{k+1}$. Let us denote it for 
short by $F$, namely $\vec{X}_{k+1}=F(\vec{X}_k)$.

Assuming the existence of a $2\pi/\nu_j$ quasiperiodic orbit corresponding to the initial condition $\vec{X}_{\infty}$ is 
equivalent to assume that $F$ admits $\vec{X}_{\infty}$ as a fixed point, i.e. $F\left(\vec{X}_{\infty}\right)=\vec{X}_{\infty}$.
Then, to prove that $\vec{X}_k$ converges to $\vec{X}_{\infty}$  we have to prove that this fixed point is indeed an attractor, 
namely that every eigenvalue of the Jacobian matrix $\Phi$ of $F$ has a modulus lower than 1.

Rewriting~\eqref{eq:sls} as

\begin{equation}
  \label{eq:sls2}
  S(0;F(\vec{X}_{k}))+L(0;F(\vec{X}_{k}))=S(0;\vec{X}_k)\, ,
\end{equation}
and then differentiate every component $S_l$ of it with respect to every component $x_i$ of $\vec{X}_k$:
\begin{equation}
  \label{eq:sls3}
  \frac{\partial S_l}{\partial x_i}=\Sigma_{j=1}^p\left(\frac{\partial F(\vec{X_k})_j}{\partial x_i}
\frac{\partial S_l}{\partial x_j}+\frac{\partial F(\vec{X_k})_j}{\partial x_i}\frac{\partial L_l}{\partial x_j}\right)
\end{equation}

\par In calling $\Sigma$ the Jacobian matrix of $\vec{S}$, we straightforwardly get from (\ref{eq:sls3})

\begin{equation}
  \label{eq:sls4}
  \Sigma(\vec{X}_k)=\Sigma\left(F(\vec{X}_k)\right)\Phi(\vec{X}_k)+\Lambda\left(F(\vec{X}_k)\right)\Phi(\vec{X}_k),
\end{equation}
what yields

\begin{eqnarray}
  \Phi(\vec{X}_k) & = & \left(\Sigma\left(F(\vec{X}_k)\right)+\Lambda\left(F(\vec{X}_k)\right)\right)^{-1}\Sigma(\vec{X}_k) \\
                  & = & \left(\Sigma\left(F(\vec{X}_k)\right)+\Lambda\left(F(\vec{X}_k)\right)\right)^{-1} \\
 & & \times\Lambda\left(F(\vec{X}_k\right)\left[\Lambda\left(F(\vec{X}_k\right)\right]^{-1}\Sigma(\vec{X}_k)
\end{eqnarray}
From the hypothesis~\eqref{eq:hypo} we have

\begin{equation}
  \lim_{\vec{X}_k\rightarrow \vec{X}_{\infty}}\Phi(\vec{X}_k)=0,
\end{equation}
so the eigenvalues of the matrix $\Phi(\vec{X}_{\infty})$ have a modulus lower than $1$.

\subsection{Discussion}

\par The caveat of this proof is that the hypothesis~\eqref{eq:hypo} usually cannot be checked a priori, so for a given dynamical system studied close to
a stable equilibrium, we cannot a priori know whether our algorithm will converge or not. Since this algorithm aims at determining the solution 
$\vec{P}=\vec{X}_{\infty}$, it is a priori unknown as well.

\par Anyway, we can prove that the hypothesis~\eqref{eq:hypo} is verified in the Hamiltonian framework, since it is in fact a weakening of the d'Alembert characteristic.
For sake of clarity, we here reduce to 1-degree of freedom systems, but the principle is the same for 
multi-dimensional systems. In this context, the hypothesis~\eqref{eq:hypo} reads

\begin{equation}
  \label{eq:hypo2}
  \lim_{x\rightarrow x_{\infty}}\frac{\partial_xS(0;x)}{\partial_xL(0;x)}=0\, .
%   \lim_{x\rightarrow x_{\infty}}\left[\Lambda(x)\right]^{-1}\Sigma(x)=0\, .
\end{equation}
Thank to the d'Alembert characteristic we have:
  \begin{equation*}
    S(0;x)\sim x_{\infty}+a|x-x_{\infty}| + \dots \quad \text{and}\quad
    L(0;x)\sim b\sqrt{|x-x_{\infty}|} + \dots \, .
  \end{equation*}
Thus $\partial_xS(0;x)$ is finite while
$\partial_xL(0;x)$ diverges, as $\sqrt{|x-x_{\infty}|}$, and
condition~\eqref{eq:hypo} is satisfied. 

More generally in the case 
  \begin{equation}
\label{eq:ab}
    S(0;x)\sim x_{\infty}+a(x-x_{\infty})^{\alpha} + \dots \quad \text{and}\quad
    L(0;x)\sim b(x-x_{\infty})^{\beta} + \dots \, ,
  \end{equation}
for some positive $\alpha$ and $\beta$, such that $\alpha > \beta$ satisfies the condition~\eqref{eq:hypo}, one can provide also 
the rate of convergence of the algorithm. In fact, let us define $\theta_n=x_n-x_{\infty}$,
then using~\eqref{eq:ab}, Eq.~\eqref{eq:sls} rewrites:
\begin{equation*}
  a\theta_{n+1}^{\alpha}+b\theta_{n+1}^{\beta}\sim a\theta_{n}^{\alpha}\, ,
\end{equation*}
that can be approximately solved for $\theta_{n+1}$, for instance using the
Lagrange inversion formula~\cite{Carletti2003}, to get
\begin{equation*}
  \theta_{n+1}\sim
  \left(\frac{a}{b}\right)^{1/\beta}\theta_n^{\alpha/\beta}+\dots \, ,
\end{equation*}
and thus providing a convergence rate $\alpha/\beta$. In the above mentioned
Hamiltonian framework this results into a quadratic convergence rate. But to have the convergence of the algorithm, we just need $\beta>\alpha$.

\section{Application of the algorithm to a problem in mathematical biology}
\label{sec:applbio}

\par In this section we want to investigate the applicability of the algorithm beyond the Hamiltonian framework. We thus chose an example from mathematical biology, a prey-predator system, derived from the Lotka-Volterra equations, periodically forced by a sinusoidal term, accounting for the seasonal changes.

\subsection{Forced prey-predator systems}

\par This model presented in Blom et al.~\cite{bbgv81} aims at representing the densities of preys $x_1$ and predators $x_2$ as a function of time; beside the standard interactions among predator and prey and the logistic growth rate for the prey, we also take into account an external periodic forcing term, that can represent the one-year-periodic density of food availability due to the seasonal alternation.  

\begin{figure}[ht]
  \centering
  \includegraphics[width=\textwidth]{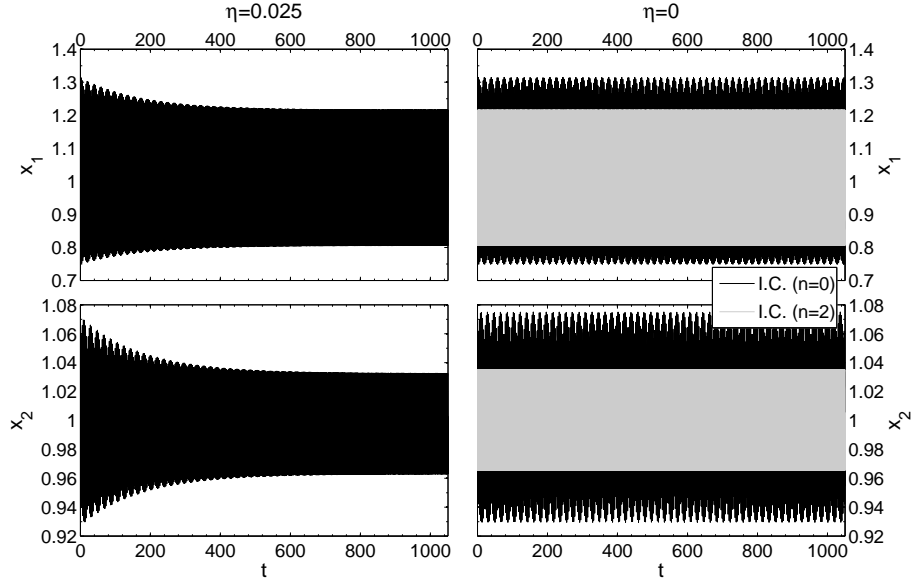}
  \caption{Numerical integration of the system of Blom et al. for $\alpha=4.539$, $\beta=1.068$ and $\gamma=0.25$, the initial conditions being $x_1=1$ and $x_2=1-\eta$ (in black), and the period $T=2.853\,74$. We can see a damping of the oscillations for $\eta=0.025$. On the right, the grey curves result from numerical integrations of the equations of the system with initial conditions obtained after 2 iterations ($n=2$). The resulting behavior is very similar to the one obtained after damping. \label{fig:x1x2}}
\end{figure}

The equations ruling the densities of the populations read:
\begin{eqnarray}
	\frac{dx_1}{dt} & = & \alpha x_1\big(1+\gamma\cos(2\pi t)-x_2-\eta x_1\big), \label{eq:dx1} \\
	\frac{dx_2}{dt} & = & \beta x_2(-1+x_1), \label{eq:dx2}
\end{eqnarray}
where $\alpha$, $\beta$, $\gamma$, $\delta$ and $\eta$ are non negative constants. When $\gamma=0$ (autonomous system) and $\eta=0$, one can straightforwardly prove the existence of $T=2\pi/\sqrt{\alpha\beta}$ periodic solutions, around the non trivial equilibrium $(x_1,x_2)=(1,1)$. If $\eta>0$ the non trivial equilibrium moves to $(x_1,x_2)=(1,1-\eta)$ and moreover these oscillations are damped; for $\gamma>0$ the system is subjected to $2\pi$-periodic forced oscillations (see Fig.~\ref{fig:x1x2} left panels). 

Using the algorithm, we can determine the initial conditions corresponding to the solution with zero amplitude (see Fig.~\ref{fig:x1x2} right panels, black curves) associated with the proper period $T$, without using any damping (i.e. with $\eta=0$).

\subsection{Numerical analysis}

\par We performed a numerical analysis of the above forced predator-prey model using the following parameters values: $\alpha=4.539$, $\beta=1.068$, $\gamma=0.25$ and $\eta=0$. 
We adopted the variable step size Bulirsh-Stoer algorithm \cite{bulirsh-stoera,bulirsh-stoerb} to numerically integrate the differential equations (\ref{eq:dx1}) and 
(\ref{eq:dx2}). The proper frequency, for the case $\gamma =0$, is given by  $\sqrt{\alpha\beta}=2.201\,74$. Results are reported in Table~\ref{tab:preypredator}.

\begin{table}[h!]
  \begin{small}
  \begin{center}
    \caption{Application of the algorithm to the forced prey-predator system, with $\alpha=4.539$, $\beta=1.068$, $\gamma=0.25$ and $\eta=0$. \label{tab:preypredator}} 
    \begin{tabular}{ >{$}l<{$} >{$}c<{$} >{$}l<{$} | >{$}c<{$} >{$}c<{$} | >{$}c<{$} >{$}c<{$} }
      \text{n} & & \text{I.C.} & \# \text{freq.} & \text{Rank} & \text{Ampl.} & \omega^{\bullet}\\
      \hline
      \hline
      \multirow{2}{*}{0} & x_1 & 1.000\,000\,000\,000\,000 & 50 & 2 & 3.831\,163\dexp{-2} & 2.206\,634\\
      & x_2 & 1.000\,000\,000\,000\,000 & 50 & 1 & 1.854\,280\dexp{-2} & 2.206\,634 \\
      \hline
      \multirow{2}{*}{1} & x_1 & 0.989\,166\,714\,745\,100 & 22 & 4 & 3.573\,335\dexp{-5} & 2.207\,483 \\
      & x_2 & 0.965\,514\,795\,157\,481 & 21 & 3 & 1.729\,063\dexp{-5} & 2.207\,483\\
      \hline
      \multirow{2}{*}{2} & x_1 & 0.989\,186\,585\,234\,344\,297\,1 & 26 & 6 & 4.508\,632\dexp{-9} & 2.207\,483 \\
      & x_2 & 0.965\,545\,142\,090\,197\,513\,7 & 26 & 6 & 2.181\,634\dexp{-9} & 2.207\,483\\
      \hline
      \multirow{2}{*}{3} & x_1 & 0.989\,186\,576\,347\,806\,470\,6 & 14 & 11 & 6.524\,090\dexp{-17} & 2.207\,483 \\
      & x_2 & 0.965\,545\,142\,191\,326\,750\,5 & 14 & 11 & 3.156\,872\dexp{-17} & 2.207\,483\\
      \hline
      \multirow{2}{*}{4} & x_1 & 0.989\,186\,576\,347\,806\,470\,2 & 14 & 11 & 6.503\,088\dexp{-17} & 2.207\,483 \\
      & x_2 & 0.965\,545\,142\,191\,326\,750\,4 & 14 & 11 & 3.146\,710\dexp{-17} & 2.207\,483\\
    \end{tabular}
  \end{center}
  \end{small}
\end{table}

We see variations of the detected proper frequency $\omega^{\bullet}$, before converging to $2.207\,4$, that is close to the predicted value of $2.201\,74$ (with $\gamma=0$). 

In right panels of Figure~\ref{fig:x1x2}, we plot (grey curves) the evolution of the variables $x_1$ and $x_2$ for the initial condition corresponding to the iteration $n=2$ (Table~\ref{tab:preypredator}). We notice that, with these initial conditions, we removed the free libration parts of the signal. The resulting behavior is very similar to the damping case ($\eta\neq 0$) but with the damping the equilibrium is shifted from $(1,1)$ to $(1,1-\eta)$.

\section{Alternative algorithms}
\label{sec:alternate}

\par In this section we briefly review other existing algorithms having goals similar to the one we present here. We do not intend to make a rigorous and exhaustive 
comparison between them and the presented algorithm because we think that these studies go beyond the scope of this paper. We thus deserve this 
for forthcoming papers. We here focus on 2 algorithms: the iterative method by Rodionov et al. (\cite{rs06,ro08,ras09,ra11}), and the fit of libration centers by Bois \& Rambaux \cite{br07}.

  \subsection{The iterative method (Rodionov et al. 2006-2011)}

\par This method has been elaborated in the framework of galactic dynamics, to find the initial conditions of N-body simulations corresponding to equilibrium states. These equilibrium states are characterized by various parameters that describe for instance the mass density of a galaxy, or the radial velocity dispersion. This algorithm consists in choosing a priori suitable initial conditions, letting the system evolve in propagating its orbit, possibly under some constraints, and use the final state of the system to build new initial conditions.

\par The authors claim that this method can be adapted to any arbitrary dynamical system, and we agree with this opinion, at least for the ones modeling physical 
phenomena. A rough adaptation of this algorithm in the context of the rotation of celestial bodies in spin-orbit resonance would consist in including a dissipative 
effect in the equations of the problem, and letting the system evolve on a long enough timescale 
(Peale \cite{p05} estimates this timescale of the order of $10^5$ years for the rotational dynamics of Mercury) for the damping to act efficiently, and then to use the computed spin angle and spin velocity to start a new simulation. A way to bypass the problem of the long damping time would be to artificially accelerate 
the damping process. Nevertheless, the dissipation changes the equilibrium position of the system (see e.g. \cite{rcwk10}). A refinement of this process could be to use a fast dissipation at the first iteration, and after a slower one. Using the presented algorithm consists in making the assumption that the effect of the dissipation on the location of the equilibrium is negligible. In fact, this effect has not been observed yet except for the Moon; detecting it for Mercury with the space missions MESSENGER and Bepi-Colombo is challenging. So, this assumption can be considered as acceptable. An advantage of the iterative method is that it does not assume the existence of a $2\pi/\nu_j$ quasiperiodic orbit.

\par In the case of the prey-predator problem (Sec.\ref{sec:applbio}), ths algorithm requires to set the parameter $\eta$ to $0$, while the iterative method would use its strictly positive value. As we have already seen, it shifts the mean equilibrium from $(1,1)$ to $(1,1-\eta)$.

  \subsection{Fit of libration centers (Bois \& Rambaux 2007)}

\par This algorithm has been elaborated in the context of a numerical computation of Mercury's equilibrium obliquity, and to the best of our knowledge has neither been used in any other study. This problem is a 2-dimensional spin-orbit problem of a rigid body, with the notable difficulty that it contains a very long period due to the regressional motion of Mercury's orbital ascending node, whose period is of the order of $250,000$ years. An efficient numerical identification of this sinusoidal term would require to propagate the equations of the system over at least $5\dexp{5}$ years (i.e. two periods), while the space missions cover a period of about 2 years, and the orbital ephemerides available cover a few thousands of years. So, trying to identify the period of regression of the node seems not to be the right way to do. Anyway, the use of unoptimized initial conditions induces 1,000-y free oscillations whose amplitude is expected to be null in the real system. So, Bois \& Rambaux propose to fit a sinusoid to the orbits obtained after numerical simulations (in which the dissipation is neglected) to identify the free oscillations, and to remove them from the initial conditions to perform a new simulation. This way of identifying free oscillations without making a complete frequency analysis of the system does not require the final solution to be $2\pi/\nu_j$ quasiperiodic (in fact, the $250,000$-y forced oscillation is seen as the slope over this timescale). So, this is just a partial decomposition of the signal, a complete one, when possible, is of course more accurate.

\section{Conclusion}
\label{sec:concl}

\par We hereby presented an algorithm able to determine the initial conditions corresponding to forced oscillating equilibria, namely removing free proper librations. 
In our restricted group of local collaborators, we are used to name it NAFFO, for Numerical Algorithm For Forced Oscillations, since this recalls its proximity to NAFF. But this is 
just for internal use, we do not claim
for the fatherhood of this algorithm and every user is of course free to give it the nickname he prefers. Given an initial condition, this algorithm iteratively 
produces better and better approximated initial conditions, by integrating the system, identifying the free and forced frequency terms and eventually removing 
the former ones. In the present paper this step has been performed using the NAFF method by Laskar.

\par Under suitable conditions, we proved the convergence of the algorithm, under the assumptions of exact numerical integration and frequency identification. We shown that in the Hamiltonian framework the required hypothesis are satisfied, whenever the d'Alembert characteristic holds true.

We provided one relevant application to benchmark the proposed algorithm, in mathematical biology. This supplements the applications already present in the literature.

\section*{Acknowledgments}
We are indebted to Julien Dufey for his help in the computation of the Lie transforms, and to Philippe Robutel and Lia Athanassoula for fruitful discussions.


\appendix

\section{The NAFF algorithm}
\label{ssec:naff}

\par Since the algorithm we present requires an accurate quasi-periodic representation of the orbit, we hereby shortly present the method we used: the NAFF (Numerical Analysis of the Fundamental Frequencies) algorithm due to J. Laskar (see e.g. \cite{l93,l05} where its convergence has been proved under suitable Diophantine hypothesis for the frequency vector). The good accuracy properties of the NAFF, allows to use it for instance to characterize  the chaotic diffusion in dynamical systems~\footnote{In this framework, it is  also known as FMA for Frequency Map Analysis.}~\cite{lfc92}, in particles accelerator~\cite{nl03} or in celestial mechanics~\cite{nv07,ldv09}.  

The algorithm aims at numerically identifying the coefficients $a_l$ and $\nu_l$ of the quasi--periodic complex signal $x(t)$ known over a finite, but large, time span $[-T;T]$, of the form~\eqref{eq:x}, thus providing an approximation
\begin{equation}
  \label{equ:naff}
  x(t) \approx
  \sum_{l=1}^Na_l^{\bullet}\exp\left(\imath f_l^{\bullet}t\right)\, ,
\end{equation}
where $\nu_l^{\bullet}$, respectively $a_l^{\bullet}$, are the numerically (as emphasized by the bullet) determined real frequencies, respectively complex coefficients, i.e. amplitudes. If the signal $x(t)$ is real, its frequency spectrum is symmetric and the complex amplitudes associated with the frequencies $\nu_l^{\bullet}$ and $-\nu_l^{\bullet}$ are complex conjugates. 

The frequencies and associated amplitudes are found with an iterative scheme. To determine the first frequency $\nu_1^{\bullet}$, one searches for the maximum of the amplitude of  
\begin{equation}
  \label{equ:philas}
  \phi(\nu)=\langle x(t),\exp(\imath\nu t)\rangle\, ,
\end{equation}
where the scalar product $\langle f(t),g(t)\rangle$ is defined by
\begin{equation}
  \label{equ:prodscal}
  \langle f(t),g(t)\rangle =\frac{1}{2T}\int_{-T}^T
  f(t)\,\overline{g(t)}\,\chi(t)\, dt, 
\end{equation}
where $\overline{g(t)}$ is the complex conjugate of $g(t)$ and where $\chi(t)$ is a weight function, i.e. a positive function verifying
\begin{equation}
  \label{equ:poids}
  \frac{1}{2T}\int_{-T}^T \chi(t)\,dt=1.
\end{equation}
Laskar advises to use
\begin{equation}
  \label{eq:chip}
  \chi(t)=\frac{2^p(p!)^2}{(2p)!}\Big(1+\cos (\pi t)\Big)^p,
\end{equation}
where $p$ is a positive integer. In practice, the algorithm is the most efficient with $p=1$ or $p=2$. We used $p=2$ in the numerical applications presented later in the paper, because it yields good results. 

Once the first periodic term $\exp(\imath\nu_1^{\bullet}t)$ is found, its complex amplitude $a_1^{\bullet}$ is obtained by orthogonal projection, and the process is started again on the remainder $f_1(t)=f(t)-a_1^{\bullet}\exp(\imath\nu_1^{\bullet}t)$. The algorithm stops when two detected frequencies are too close to each other, because this alters their determinations, or when the number of detected terms reaches a maximum set by the user. When the difference between two frequencies is larger than twice the frequency associated with the length of the total time interval, the determination of each fundamental frequency is not perturbed by the other ones. Once the frequencies have been determined, it is also possible to refine their determination iteratively, numerically \cite{c98} or analytically \cite{sn97}.

\end{document}